\documentclass[graybox]{svmult}

\usepackage{mathptmx}       
\usepackage{amsmath}
\usepackage{amssymb}
\usepackage{helvet}         
\usepackage{courier}        
\usepackage{type1cm}        
\usepackage{makeidx}         
\usepackage{graphicx}        
\usepackage{multicol}        
\usepackage[bottom]{footmisc}

\usepackage{natbib}

\newcommand{\cG}{{G}}
\newcommand{\cA}{{A}}
\newcommand{\cN}{{N}}
\newcommand{\cNini}{\cN^{-}_i}
\newcommand{\cNouti}{\cN^{+}_i}

\newcommand{\cM}{{ M}}
\newcommand{\cT}{{ T}}
\newcommand{\cK}{{ K}}

\newcommand{\xak}{{x_{a}^k}}
\newcommand{\xijk}{{x_{ij}^k}}
\newcommand{\xjik}{{x_{ji}^k}}

\newcommand{\tki}{{t_{ki}}}
\newcommand{\tkj}{{t_{kj}}}
\newcommand{\tij}{{t_{ij}}}

\newcommand{\poly}{\hbox{$ P^{\hbox{\tiny ND}}$}}
\newcommand {\beq}{\[}\newcommand {\eeq}{\]}
\newcommand {\beqn}{\begin{equation}}\newcommand {\eeqn}{\end{equation}}
\newcommand {\beqan}{\begin{eqnarray}}\newcommand {\eeqan}{\end{eqnarray}}
\newcommand {\beqa}{\begin{eqnarray*}}\newcommand {\eeqa}{\end{eqnarray*}}

\newcommand{\R}{\ensuremath{\mathbb{R}}}
\newcommand{\cred}{\color{red}}
\newcommand{\Z}{\ensuremath{\mathbb{Z}}}

\newcommand{\be}[1]{\begin{equation}\label{#1}}
\newcommand{\ee}{\end{equation}}

\newcommand{\conv}{\ensuremath{\text{conv}}}

\makeindex             

\newcommand{\ignore}[1]{}
\def\FSS{\ensuremath{{F}_{SS}} }
\def\FSM{\ensuremath{{F}_{SM}} }
\def\FMS{\ensuremath{{F}_{MS}} }
\def\FMM{\ensuremath{{F}_{MM}} }
\newcommand{\RR}{\R}
\newcommand{\ZZ}{\Z}
\newcommand{\st}{\ensuremath{\text{s.t. }}}
\newcommand{\NP}{\ensuremath{\mathcal{NP}}}

\newcommand{\floor}[1]{\ensuremath{\lfloor #1 \rfloor}}
\newcommand{\ceil}[1]{\ensuremath{\lceil#1\rceil}}
\newcommand{\xceil}[1]{\ensuremath{\left\lceil#1\right\rceil}}
\newcommand{\FLOOR}[1]{\ensuremath{\lfloor #1 \rfloor}}
\newcommand{\CEIL}[1]{\ensuremath{\lceil#1\rceil}}
\usepackage[colorinlistoftodos,prependcaption,textsize=small]{todonotes}

\newcommand{\update}[1]{{\color{black}#1 }}

\newcommand{\FU}{\ensuremath{F_U}}
\newcommand{\FUP}{\ensuremath{F_U(\FIXY, K_0, K_1)}}

\newcommand{\FUr}{\ensuremath{F_{Ur}}}
\newcommand{\FS}{\ensuremath{F_S}}
\newcommand{\FB}{\FU}
\newcommand{\FC}{\FS}

\newcommand{\FIXETA}{\CEIL{\bar y}}
\newcommand{\FIXY}{\ensuremath{\nu}}
\newcommand{\frcpnt} {\ensuremath{( \bar x, \bar y )}}

\newcommand{\rc}{residual capacity }
\newcommand{\sm}{\setminus}
\newcommand{\emptybox}[1]{\begingroup	\setlength{\fboxsep}{-\fboxrule}\noindent\mbox{\rule{0pt}{#1}}	\endgroup}

\begin{document}

\title*{Multi-Commodity Multi-Facility Network Design}

\author{Alper Atamt\"urk and Oktay G\"unl\"uk}
\institute{
	Alper Atamt\"urk  \at Department of Industrial Engineering and Operations Research, University of 	California, Berkeley, CA 94720 USA. \email{atamturk@berkeley.edu}
\and 
Oktay G\"unl\"uk \at School of Information Engineering and Operations Research, Cornell University, Ithaca, NY 14853. \email{oktay.gunluk@cornell.edu}}

%
%
\maketitle

\abstract{
		We consider multi-commodity network design models, where capacity can be added to the arcs of the network using multiples of facilities that may have different capacities. This class of mixed-integer optimization models appears frequently in telecommunication network capacity expansion problems, train scheduling with multiple locomotive options, supply chain and service network design problems. 
	Valid inequalities used as cutting planes in branch-and-bound algorithms have been
	instrumental in solving large-scale instances. We review
	the progress that has been done in polyhedral investigations in this area by emphasizing three fundamental techniques. These are the 
	metric inequalities for projecting out continuous flow variables, mixed-integer rounding from appropriate base relaxations, and shrinking the network to a small $k$-node graph.	 
	The basic inequalities derived from arc-set, cut-set and partition relaxations of the network
	are also extensively utilized with certain modifications in robust and survivable network design problems.
}

\section{Introduction}
Here we consider multi-commodity network design models, where capacity can be added to arcs of the network using integer multiples of facilities, possibly with varying capacities. This class of models appear frequently in telecommunication network capacity expansion problems
train scheduling with multiple locomotive options,
supply chain and service network design problems. 
In the single-facility network design problem, one installs multiples of only
a single type of facility on the arcs of the network. Routing vehicles with identical
capacity in a logistics network and installing a communication network with only one
cable type are examples of the single-facility network design problem. In the multi-facility
problem, one may install different types of facilities with varying capacities,
such as fiberoptic cables with varying bandwidths, production lines or machines with
different rates, or a fleet of heterogeneous vehicles with varying capacities.
The optimization problem seeks to decide how many facilities of each type to install 
on the network so as to meet the demand for each commodity at the least cost.
We present the precise problem description and the associated formulation in the next section.

Different versions of the problem are obtained depending on how the flow is routed in the network.
In the {\em unsplittable} flow version, only a single path is allowed to route the flow from its source to its destination, which requires integer variables to model its route. 
This is the case, for instance, in telecommunication networks using multiprotocol label switching (MPLS) technique,
production and distribution with single sourcing,
and express package delivery. 
The {\em splittable} case, which assumes that flow can be routed using  multiple directed paths,
is obviously a relaxation of the unsplittable case; therefore, valid inequalities for the splittable case are also valid for the unsplittable case.

In addition, the capacity created by the facilities can be {\em directed},  {\em bidirected}, or {\em undirected}.
In the bidirected case, if a certain facility is installed on an arc, then the same facility also needs to be installed on the reverse arc. In the undirected case, the total flow on an arc and its reverse arc is limited by the capacity of the undirected edge associated with the two arcs. Here we consider the directed case where the total flow on an arc is limited by the total (directed) capacity of the arc. 
In a recent paper, the authors of this chapter describe how to generate valid inequalities  for the bidirected and undirected network design problems using valid inequalities for the directed case.

In this chapter we review strong  valid inequalities for the multi-commodity, multi-facility network design problem. Throughout, we emphasize three fundamental techniques that have been effective in deriving strong inequalities for network design problems. These are the metric inequalities for projecting out the continuous flow variables, mixed-integer rounding from appropriate base relaxations, and shrinking the network to a small $k$-node graph. The valid inequalities for the network design problem are obtained by applying these techniques to different relaxations of the problem, namely, arc-set, cut-set and partition relaxations.
The basic inequalities derived from these relaxations
are also utilized, with certain adaptations, in robust and survivable network design problems.

This chapter is organized as follows. In the next section, we introduce the notation used in the chapter and 
give a formal definition of the multi-commodity multi-facility design problem. 
Section~\ref{sec:prelim} reviews the metric inequalities for projecting out the multi-commodity flow variables, the
mixed-integer rounding technique, as well as a simplification obtained by shrinking the network for deriving valid inequalities from a smaller network. These 
techniques play a central role in deriving strong valid inequalities for the network design problem. 
Section~\ref{sec:arcset} reviews the inequalities from single arc capacity constraints for the splittable as well as the unsplittable cases.
Section~\ref{sec:cutset} reviews the valid inequalities from two-partitions for single as well as multi-facility cases.
Section~\ref{sec:part} generalizes the inequalities in the previous section to a higher number of partitions.
We conclude with a few final remarks in Section~\ref{sec:conc}.

\section{Problem Formulation}
\label{key}


Let $\cG =(\cN,\cA)$ be a directed graph (network), with node set $\cN$ and arc set $\cA \subseteq \cN\times \cN$.
Each arc $a\in\cA$  has a given existing capacity $\bar c_a\ge0$ and the network design problem consists of installing additional capacity on the arcs of the network using an (integral) combination of a given set of capacity types.
The objective of the problem is to minimize the sum of the flow routing cost and 
the capacity expansion cost.
Throughout we assume that the data is rational.

The demand data for the problem is given by matrix $\cT=\{\tij\}$, where $\tij\ge0$ is
the amount of (directed) traffic that should be routed from node $i\in \cN$ to $j\in \cN$. 
Using matrix $\cT$, we can define a collection of commodities, each of which has a certain supply and demand at each node of the network. 
As discussed in Chapter 2, there are different ways to formulate the same network design problem by changing what is meant by a commodity. 
Using a minimal vertex cover of the so-called demand graph, one can  obtain the smallest number of commodities required to formulate a given problem instance correctly.
Computationally, formulations with fewer commodities may be more desirable as their continuous relaxations  solve faster. 
For the sake of concreteness,  we will use the aggregated commodity description but we note that most of the discussion below does not depend on how the commodities are defined. 
Let $\cK\subseteq\cN$ denote the collection of nodes with positive supply, i.e., 
$$\cK= \bigg \{i\in\cN : \sum_{j \in \cN}\tij>0 \bigg \} \cdot$$
We use $ w_i^k$ to denote the net demand of node $i\in\cN$ for commodity $k\in\cK$.
More precisely, let  $ w_i^k=\tki$ for $i\not=k$ and  $w_k^k = -\sum_{j \in \cN} \tkj$ for $k\in\cK$.
Note that each node $k\in \cK$ is the unique supplier of commodity $k$ and flow of each commodity in the network needs to be disaggregated to obtain  an individual routing for  origin-destination pairs.
We also note that the aggregated commodity description can only be used if the flow routing cost  does not depend on both the origin and the destination of the traffic that needs to be routed.

New capacity can be installed in the network using integer multiples of facilities $\cM$,
where a single unit of facility $m\in\cM$ provides capacity $c_m$. 
Without loss of generality, we  assume that $c_m\in\Z$ for all $m\in\cM$ and $c_1<c_2<\ldots<c_{|M|}$.
In this setting the network design problem involves installing enough additional capacity on the arcs of the network so that traffic can be routed simultaneously without violating arc capacities.
For $i\in \cN$, let
\begin{equation} \nonumber
\cNouti = \{ j \in \cN : (i, j) \in \cA \} \text{ and } 
\cNini = \{ j \in \cN : (j, i) \in \cA \}
\end{equation}
denote the neighbors of node $i\in\cN$. 
Let integer  variables $ y_{ma}\ge0$ indicate the number of facilities of type $ m\in\cM$ installed on arc  $a \in \cA$ and continuous variables $\xak\ge0$ denote the amount of flow of commodity $k\in\cK$ routed on arc  $a \in \cA$.  Using this notation, the following constraints define the feasible region of the multi-commodity multi-facility network design problem:
\begin{align}
\label{eq:balance}
\sum_{j \in \cNouti} \xijk ~~- \sum_{j \in \cNini} \xjik~~~~  = ~~w_i^k, && {\rm \ } k \in \cK, \  i\in \cN, 
\\
\label{eq:cap}
\sum_{k \in \cK } x_{ij}^k\ - \sum_{m \in \cM } \ c_m  y_{mij}  ~\leq~~ \bar c_{ij}, &&{ ~~~\ } (i,j) \in\cA.   
\end{align}
Then the network design problem is stated as: 
\beqn \hskip -3cm \text{(NDP)}  \hskip 2cm  \text{min} \Big\{dy+fx\::\: (x,y)\in\poly\Big\}, \nonumber
\eeqn
where $d$ and $f$ are cost vectors of appropriate size and
$$\poly =  \conv \Big\{   (x,y) \in     \R^{\cA \times \cK}_+ \times \Z^{\cA\times\cM}_+:~ \eqref{eq:balance} \text{ and } \eqref{eq:cap} \Big \} \cdot$$

\begin{figure}[h]\centering
\begin{tikzpicture}[scale = 0.5]
\draw[thick,->] (0,0) -- (18,0) node[anchor=north west] {$z$};
\draw[thick,->] (0,0) -- (0,9) node[anchor=south east] {$h(z)$};

\foreach \j in {0,1} {	
	\foreach \i in {0,...,2} {
		\draw (7*\j+\i+0.07,\j*3.75+\i+1) -- (7*\j+1+\i,\j*3.75+\i+1);
		\draw (7*\j+\i,\j*3.75+\i+1) circle (0.07);
		\fill (7*\j+1+\i,\j*3.75+\i+1) circle (0.07);
	}
	\draw(7*\j+3+0.07,\j*3.75+3.75)--(7*\j+7,\j*3.75+3.75);
	\draw(7*\j+3,\j*3.75+3.75) circle (0.07);
	\fill(7*\j+7,\j*3.75+3.75) circle (0.07);
}
\fill(0,0) circle (0.07);

\draw (1,0.1)--(1,-0.1); \draw (1,-0.2)node[anchor=north] {$ c_1$};
\draw (2,0.1)--(2,-0.1) node[anchor=north] {$ 2c_1$};
\draw (3,0.1)--(3,-0.1) node[anchor=north] {$ 3c_1$};
\draw (7,0.1)--(7,-0.1); \draw (7,-0.2) node[anchor=north] {$c_2$};
\draw (8,0.1)--(8,-0.1) ; \draw (8.5,-0.1) node[anchor=north] {$c_2+c_1$};

\draw (0.1,1)--(-0.1,1) node[anchor=east] {$  d_{1,a}$};
\draw (0.1,2)--(-0.1,2) node[anchor=east] {$ 2d_{1,a}$};
\draw (0.1,3)--(-0.1,3) node[anchor=east] {$ 3d_{1,a}$};
\draw (0.1,3.75)--(-0.1,3.75) node[anchor=east] {$ d_{2,a}$};
\draw (0.1,4.75)--(-0.1,4.75) node[anchor=east] {$ d_{1,a}+d_{2,a}$};

\end{tikzpicture}
\caption{The piecewise linear capacity installation cost function $h(z)$.} 	\label{fig:economiezofscale}
\end{figure}
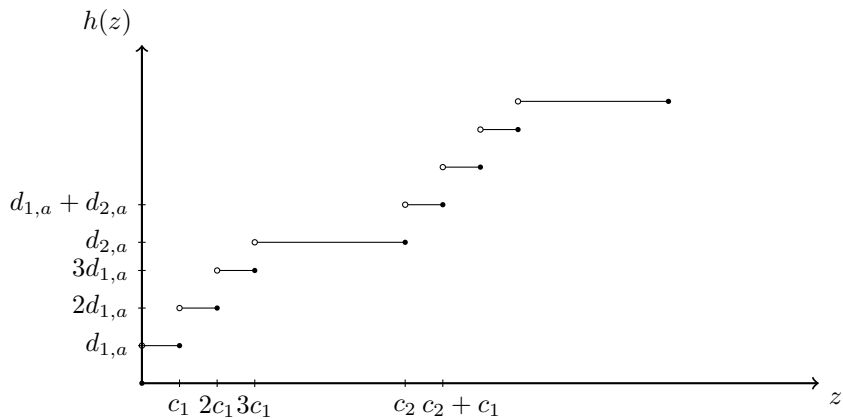

As a concrete example with two facilities, consider a given arc $a\in \cA$.
The total capacity $c_1 y_{1a}+ c_2  y_{2a}$ given by the integer variables $y_{1a}$ and $y_{2a}$ has cost  $ d_{1a}  y_{1a}+ d_{2a}  y_{2a}$.
Assuming economies of scale, let $d_{1a}/c_1 > d_{2a}/c_2$ and remember that $c_1<c_2$.
In this case, we can write the cost function $f(z)$ required to generate $z$ units of capacity (at the least cost) as: 
$$h(z)= \floor{z/c_2}d_{2a}+\min\{d_{2a},\ceil{(z- \floor{z/c_2}c_2)/c_1}d_{1a}\}$$
which is a piecewise linear function.
Figure~\ref{fig:economiezofscale} illustrates an example with $3d_{1a}<d_{2a}<4d_{1a}$.

We also note that, when $f=0$ it is  possible to project out the multi-commodity flow variables from \poly\ to
obtain a formulation in the space of only the discrete capacity variables. 
This capacity formulation requires an exponential number of
constraints and is discussed in Section \ref{sec:prelim-metric}.

\ignore{
	\begin{alignat*}{2}
	\min \sum_{a \in A, k \in K} e^k_a f^k_a + \sum_{a \in A, t \in T} & d^t_a x^t_a \\
	\sum_{a \in \delta^+(v)} f_{a}^k - \sum_{a \in \delta^-(v)} f_{a}^k  & =  b^k_v, &  \  & \forall v \in V, \ \forall k \in K,\\
	\sum_{k \in K} f^k_a & \leq  \sum_{t \in T} c^t x^t_a,  & \  & \forall a \in A, \\
	f_{a}^k  \in  \RR_+, \ \forall a \in A, \ \forall k & \in K;    \ x_{a}^t \in  \ZZ_+, & &  \forall a \in A, \ \forall t \in T.
	\end{alignat*}
}

\newcommand {\Q}{{Q}}
\def\iref#1{(\ref{#1})}\def\iiref#1{inequality~\iref{#1}}
\def\Iiref#1{Inequality~\iref{#1}}
\newcommand {\bhat}{\hat b}\newcommand {\chat}{\hat c}
\newcommand {\bbar}{\bar b}\newcommand {\cbar}{\bar c}
\newcommand {\bat}{\begin{tabular}}\newcommand {\eat}{\end{tabular}}


\section{Preliminaries}\label{sec:prelim}
In this section we discuss three fundamental approaches 
that are useful in generating strong cutting planes for \poly.
We start with the metric inequalities, which give a generalization of the well-known max-flow min-cut theorem to multi-commodity flows. 
We then  discuss how valid inequalities can be generated by shrinking the network to 
one with a few nodes to obtain inequalities from simpler sets. 
Finally, we describe the mixed-integer rounding procedure, 
which is an effective method to produce valid inequalities for general mixed-integer sets.
\subsection{Metric Inequalities}\label{sec:prelim-metric}
Metric inequalities are  introduced by  \cite{I:metric}, and  \cite{OK:metric} as a generalization of the max-flow min-cut duality to  multi-commodity flows. 
Consider the following capacitated multi-commodity flow set 
\beqa
F= \bigg\{ x\in \R_+^{\cA\times\cK} \::\:  
\eqref{eq:balance} \text{ and } \sum_{k\in \cK} x_a^k   &\le&  c_{a},   ~~~~{\rm \ } a \in \cA  \:\bigg\},
\eeqa
where  $ c\in\R^{\cA}_+$ denotes the arc capacities.
By Farkas' Lemma,  the set $F $ is non-empty if and only if the following {\em metric inequalities} 
\beqan  \sum_{a\in \cA} c_a v_a &\ge& \sum_{ i\in \cN } \sum_{ k\in \cK }  w_i^ku_{i}^k \label{meta}\eeqan
hold for all $ (v,u) \in D$, where
\beqa  D= \bigg\{ (v,u) \in \R_+^{\cA}\times \R^{\cN\times\cK}\::\:
v_{ij}	&\ge&  u_{j}^k - u_{i}^k, \:\:\: k\in K, \: (i,j) \in\cA ,\  u_{k}^k=0\:\bigg\} \cdot
\eeqa
In other words, the arc capacity vector $ c$ can accommodate a feasible routing of the commodities if and only if it satisfies all metric inequalities generated by the non-empty cone $D$.
Note that for any fixed $\bar v\in \R_+^{A}$, a point $(\bar v, u(\bar v))\in D$ maximizes the right-hand-side of (\ref{meta}) when $u(\bar v)$ corresponds to the shortest path lengths using $\bar v$ as edge weights (hence the name ``metric inequality'').
Therefore, it suffices to consider metric inequalities where the vector $u \in \R^{\cN\times\cK}$ satisfies this property.

When there is only a single commodity, the max-flow min-cut theorem gives a nice characterization of 
the important extreme rays of the cone $D$. 
More precisely, in this case it suffices to consider vectors $ v\in \{0,1\}^{A}$, where $v_{ij}=1$ if and only if $i\in S$ and $j\not\in S$ for some $S\subset N$ that contains the source of the commodity but not all the sinks.

Now consider the set \poly\ and a given capacity vector $y\in\R^{\cA\times\cM}$ together with the existing arc capacities $\bar c$ and demand $w$.
The metric inequality generated by $(u,v)\in D$ becomes
\beqn  \sum_{a\in \cA}(\bar c_a +\sum_{m \in \cM } \ c_m  y_{ma} )v_a \ge \sum_{ i\in \cN } \sum_{ k\in \cK }  w_i^ku_{i}^k\label{eq:metax}.\eeqn
As it is possible to check if there is a violated metric inequality in polynomial time (by solving a linear program), one can project out the flow variables from \poly\ and obtain a ``capacity" formulation in the space of the capacity variables only. 
Clearly this approach can be applicable only if there is no flow routing cost, i.e., $f=0$ in problem (NDP).
We also note that as inequalities \eqref{eq:metax} do not have flow variables, they only depend on the demand matrix and not on what commodity definition is used for the flow variables.
Consequently, the right-hand-side of inequality \eqref{eq:metax} reduces to $\sum_{ i\in \cN } \sum_{ k\in \cN }  t_{ki}u_{i}^k$.

The {\em integral metric inequalities} are obtained by rounding up the right-hand-side of metric inequalities associated with integral vectors $(u,v)\in D$,
\beq  \sum_{a\in \cA}\sum_{m \in \cM } \ c_m  y_{ma} v_a \ge 
\xceil{\sum_{ i\in \cN } \sum_{ k\in \cK }  w_i^ku_{i}^k-\sum_{a\in \cA}\bar c_a v_a}
\label{eq:imetax}.\eeq

The basic cut-set inequalities 
discussed in Section~\ref{sec:cutset} are special cases of the integral metric inequalities.
While the metric inequalities
different from the cut-set inequalities (Section~\ref{sec:cutset}) can be useful in strengthening the convex relaxations, their separation requires more computational effort.

\subsection{Node Partition Inequalities}\label{sec:prelim-part}
Consider a partition of the node set of the directed graph $\cG =(\cN,\cA)$ into $p<|N|$ disjoint subsets: $N=\cup_{t=1}^pN_t$.
By shrinking these node subsets into singleton nodes, one obtains a simplified directed graph $ \tilde\cG =( \tilde\cN, \tilde\cA)$ with  $p$ nodes  and up to $p(p-1)$ arcs. 
In this new graph, there is an arc $(i,j)\in \tilde\cA$ from node $i\in \tilde\cN$ to node $j\in  \tilde\cN$ if and only if the original graph has at least one  arc $(u,v)\in\cA$ from some node  $u\in N_i$ to a node in $v\in N_j$. 
The existing capacity $\bar c_{ij}$ on arc $(i,j)\in  \tilde\cA$ in the new network equals the sum of the existing capacities of all the arcs from the nodes in $N_i$ to the nodes in $N_j$; in other words, $ \tilde c_{ij}=\sum_{u\in N_i}\sum_{v \in N_j} \bar c_{uv}$. 

Finally, setting the demand of $j\in\tilde\cN$ for commodity $k\in K$ to the total net demand of all nodes in $N_j$ for commodity $k$ in the original problem leads to a smaller network design problem with $p$ nodes. In other words, $ \tilde w_i^k=\sum_{v\in N_i}w_v^k$ for all $k\in K$ and $i\in \tilde\cN$. 
Note that one can reduce the number of commodities in the new problem by aggregating the ones with the same source node but in order to keep the notation simple, we do not discuss it here.

Now consider a feasible (integral) solution $( x, y)$ to the original network design problem defined on $\cG =(\cN,\cA)$ with existing capacity vector $\bar c$ and commodity demands $w$. 
Aggregating the flow and capacity variables as described above, it is easy to see that the resulting flow and the capacity vector $(\tilde x,\tilde y)$ gives a feasible solution to the simplified $p$-node  problem defined on $\tilde\cG =(\tilde\cN,\tilde\cA)$ with existing capacity vector $\tilde c$ and commodity demands $\tilde w$. 
This observation implies that valid inequalities for the simplified problem on $\tilde G$ can be translated to valid inequalities for the original problem on $G$ in the following way. 
Let
\beqn\sum_{ k\in K^{\emptybox{.13cm}}}~\sum_{ (i,j)\in \tilde A}~ \tilde\alpha_{ij}^k x_{ij}^{k}+\sum_{m \in \cM^{\emptybox{.13cm}} }\sum_{ (i,j)\in \tilde A}\tilde\beta_{mij} y_{mij}~\ge ~\gamma\label{ineq_small}\eeqn
be a given valid inequality for the simplified problem on $\tilde G$.
Then the following inequality is valid for the original problem on $G$
\beqn\sum_{ k\in K^{\emptybox{.13cm}}}~\sum_{ (u,v)\in A}~\alpha_{uv}^k x_{uv}^k+\sum_{m \in \cM^{\emptybox{.13cm}} }\sum_{ (u,v)\in A}\beta_{muv} y_{muv}~\ge~ \gamma,\label{ineq_big}\eeqn
where for  any $k\in K$, $m\in M$ and $(u,v)\in A$ with $u\in N_i$ and $v\in N_j$, we set
$$
\alpha_{uv}^k=	\begin{cases}	0, & \text{if}\ i=j \\	\tilde\alpha_{ij}, & \text{otherwise}\end{cases}~~~~~~~~~~~~
\beta_{muv}=	\begin{cases}	0, & \text{if}\ i=j \\	\tilde\beta_{mij}, & \text{otherwise.}\end{cases}
$$

\ignore{
	\cite{A:Hamid:2015} show that if inequality \iref{ineq_small} is facet-defining for the network design problem on $\tilde\cG$ with $\tilde c$ and $\tilde w$, then inequality \iref{ineq_big} is facet-defining for \poly\
	provided that $\tilde\alpha=0$,~ $\gamma> 0$, and node sets $N_1,\ldots,N_p$ induce connected components of  $G$.
	In addition, \cite{RKOW:cutset} show that the same  result holds without the assumption $\tilde\alpha=0$  when $p=2$ and $|M|=1$.
}


\subsection{MIR Inequalities} \label{sec:mir}

Many valid inequalities that have found use in practice for mixed-integer optimization problems 
are based on the mixed-integer rounding (MIR) procedure of  \cite{NW:IPbook}.
\cite{W:IPbook} illustrates the basic mixed-integer rounding on the following 
two variable mixed-integer set 
$$\Q= \Big\{(x,y)  \in \R \times  \Z\::\: x  +y  ~\ge~ b,~~x  \ge0\Big\},$$ 
and shows that the {\em basic mixed-integer rounding} inequality
\beqn x + r y \ge r \ceil{b}, \label{basic-mir} \eeqn
where $r = b - \floor{b}$, is valid and facet-defining for $\Q$.
Observe that if $b$ is integer valued, inequality \eqref{basic-mir} reduces to $x \ge 0$. Otherwise,
the inequality goes through feasible points $(0,\ceil b)$ and $(r, \floor b)$, cutting off the fractional vertex $(0,b)$. This basic principle can be applied to more general mixed-integer sets defined by a single {\em base} inequality as follows.
Let 
\beq P= \Big\{x\in \R^{n},~ y\in \Z^{l}\::\:ax+cy~\ge~ b,~~x,~y\ge0\Big\} \eeq
where $a\in\R^n$ and $c\in\R^l$. Letting $r_j$ denote $c_j-\floor {c_j}$ for $j=1, \ldots, l$, we can rewrite the base inequality as
\beq 	\sum_{a_j<0} a_jx_j +\sum_{a_j>0} a_jx_j
+ \sum_{r_j< r} r_jy_j+\sum_{r_j\ge r }  r_jy_j
+\sum_{j=1}^l \floor{ c_j} y_j~\ge~ b	\eeq
and relax it by dropping the first term, which is non-positive, and increasing the coefficients of the  fourth term, which are non-negative, to obtain the valid inequality
\beq\Big(	\sum_{a_j>0} a_jx_j + \sum_{r_j< r}r_jy_j\Big)
+\Big(\sum_{r_j\ge r}y_j+\sum_{j=1}^l \floor{c_j} y_j\Big) \ge~ b.	\eeq

As the first two sums above are non-negative and the last two sums are integer valued, treating the first two as the nonnegative continuous variable in the set $Q$ and the second two as the integer variable as in $Q$, we obtain the \textit{MIR inequality}
\beqn	\sum_{a_j>0} a_jx_j 	+ \sum_{r_j<r} r_jy_j
+ r \Big(\sum_{r_j\ge r}y_j
+\sum_{j=1}^l \floor{c_j} y_j\Big)~\ge~ r \ceil{ b}\label{eq:mir}	\eeqn
for $P$.
This MIR inequality is generated from the base inequality $ax+cy\ge b$.
Notice that, given a mixed-integer set, any valid inequality for it can be used as a base inequality to define a relaxation of the original set. Consequently, any implied inequality leads to an MIR cut.
Gomory mixed-integer cuts, for example, can be seen as MIR cuts generated from base inequalities obtained from the simplex tableau.


\section{Valid Inequalities From Arc Sets} 

\label{sec:arcset}

In this section we review the strong valid inequalities obtained  from single-arc capacity relaxations of the multi-commodity network design problem.
For simplicity of the presentation, we first focus on the single-facility case. 
We consider both the \textit{splittable-flow arc set:} 
\[
\FS \, = \Big \{ (x,y) \in [0,1]^K \times \ZZ : \ \sum_{i \in K} a_i x_i \leq a_0 + y \Big \} 
\]
and the \textit{unsplittable-flow arc set:}
\[
\FU \, = \Big \{ (x,y) \in \{0,1\}^K \times \ZZ : \ \sum_{i \in K} a_i x_i \leq a_0 + y \Big \} \cdot
\]
In Section~\ref{sec:arcset-mf} we consider the multi-facility generalizations of these sets.

The set  $F_U$ arises  in unsplittable multicommodity problems where flow between each source-sink pair needs to be routed on a single path.  In these problems,  the disaggregated commodity definition is used and the set $K$ contains all node pairs with positive demand.
In the formulation,  a binary flow variable $x_{a}^k$ is used for each commodity--arc pair $(k,a)$ that takes on a value of 1 if and only if the commodity is routed through the arc.
Consequently,  for each arc of the network this formulation has a capacity constraint of the form
\begin{align}
\hskip 4mm
\label{eq:arcset}
\sum_{k \in K} d^k x_{a}^k \leq \bar c_{a} + cy_a,
\end{align}
where  $d^k > 0$ is the demand of commodity $k\in K$,
$\bar c_{a} \ge 0$ is the existing capacity, and $c > 0$ is the unit capacity to install.
One arrives at $F_U$ by dividing \eqref{eq:arcset} by $c$.

Similarly, one arrives at $F_S$ by  redefining the flow variables associated with an arc and a commodity as the fraction of the total supply of that commodity traveling on that arc. In this case flow variables take values in $[0,1]$ and capacity constraint \iref{eq:cap} takes the form \iref{eq:arcset}. Again, dividing \eqref{eq:arcset} by $c$ gives the set $F_S$.

Without loss of generality, we assume that $a_i > 0$ for all $i \in K$, since 
if $a_i < 0$, $x_i$ can be complemented and
if $a_i=0$, $x_i$ can be dropped. 

\subsection{Splittable-flow Arc Set} \label{sec:fs}

In this section we review the valid inequalities for the splittable flow arc set \FS.
For $S \subseteq K$, by complementing the continuous variables $x_i, \ i \in S$, we can restate the arc capacity inequality as
\begin{align} \label{eq:arc-base} 
\sum_{i \in S} a_i (1-x_i) - \sum_{i \in K \sm S} a_i x_i + y \geq a(S)-a_0, 
\end{align}
where $a(S)$ stands for $\sum_{i\in S}a_i$.
Relaxing the inequality by dropping $x_i, \ i \in K \sm S$ and applying the MIR inequality, we obtain
the {\it residual capacity inequality} 
\begin{align}
\hskip 4mm
\label{ineq:rc}
\sum_{i \in S} a_i (1-x_i) \geq  r(\eta - y),
\end{align}
where $\eta = \CEIL{a(S) - a_0}$ and $r = a(S)-a_0 - \FLOOR{a(S)-a_0}$.
The residual capacity
inequalities together with the inequality $\sum_{i \in K} a_i x_i \leq a_0 + y$ and variable bounds
are sufficient to describe $conv(\FC)$.

\begin{example}
	\label{sec:ex-split-arc}
	Consider the splittable arc set
	$$\FS = \bigg \{(x,y) \in [0,1]^5 \times \ZZ: \frac{1}{3}x_1+\frac{2}{3}x_2+\frac{2}{3}x_3 \leq y \bigg \} \cdot$$
	The non-dominated arc residual capacity inequalities for \FS \ with $r >0$  are 
	\[
	\begin{array}{l|c|l}
	S \ \ & \ \ r \ \ \ & \text{ Inequalities } \\ \hline
	\{1\} & 1/3 & \  \ x_1 \le y \\ 
	\{2\} & 2/3 &\  \ x_2 \le y \\
	\{3\} & 2/3 &\  \ x_3 \le y \\
	\{2,3\} &  1/3 &\  \ 2x_2 + 2x_3 \le 2 + y \\
	\{1,2,3\} & 2/3 & \  \ x_1 + 2x_2 + 2x_3 \le 1 + 2y
	\end{array}
	\]
	
	\ignore{
		$x_1 \le y$, \ 	$x_2 \le y$, $x_3 \le y$, $2x_2 + 2x_3 \le 2 + y$, and 
		$x_1 + 2x_2 + 2x_3 \le 1 + 2y$. 
		
		$\frac{1}{3} (1-x_1) \geq \frac{1}{3} (1-y)$, 
		$\frac{2}{3} (1-x_2) \geq \frac{2}{3} (1-y)$, 
		$\frac{2}{3} (1-x_3) \geq \frac{2}{3} (1-y)$, 
		$\frac{1}{3} (1-x_1) + \frac{2}{3} (1-x_2) \geq 0 (1-y)$, 	
		$\frac{1}{3} (1-x_1) + \frac{2}{3} (1-x_3) \geq 0 (1-y)$, 
		$\frac{2}{3} (1-x_2) + \frac{2}{3} (1-x_3) \geq \frac{1}{3} (2-y)$, 
		$\frac{1}{3} (1-x_1) + \frac{2}{3} (1-x_2) + \frac{2}{3} (1-x_3) \geq \frac{2}{3} (2-y)$, 
	}	
\end{example}

Given a fractional point $\frcpnt$, a violated 
residual capacity inequality, if it exists, can be found by
the following simple separation procedure:
Let $T = \{ i \in K: \bar x_i > \bar y - \FLOOR{\bar y}\}$. 
If $a_0 + \FLOOR{\bar y} < a(T) < a_0 + \FIXETA$ and 
$\sum_{i \in T}  a_i (1-\bar x_i-\FIXETA +\bar y) + (\FIXETA - \bar y) (a_0 + \FLOOR{\bar y}) < 0$, 
then the inequality 
$\sum_{i \in T} a_i (1-x_i) \geq r (\eta - y)$ is violated by \frcpnt.
Otherwise, there exists no \rc inequality violated by \frcpnt.
Clearly, this procedure can be performed in linear time.

\subsection{Unsplittable-flow Arc Set}
\label{sec:fu}

In this section we review the valid inequalities for the unsplittable flow arc set \FU. First, consider 
the related set
\[
\FUr \, = \Big \{ (x,y) \in \{0,1\}^K \times \ZZ : \ \sum_{i \in K} r_i x_i \leq r_0 + y \Big \},
\]
where $r_i = a_i - \floor{a_i}$, $i \in K \cup \{0\}$.  
There is a one-to-one relationship between the facets of 
$\conv(\FU)$ and $\conv(\FUr)$. In particular, inequality  $\sum_{i \in K} \pi_i x_i \le \pi_0 + y$ defines a facet for \conv(\FU) \ 
if and only if inequality $\sum_{i \in K} (\pi_i - \floor{a_i}) x_i \le \pi_0 - \floor{a_0}+ y$ defines a facet for \conv(\FUr). 
Therefore, we may assume, without loss of generality, that 
$0 < a_i < 1$ for all $i \in K$ and $0 < a_0 < 1$.

\subsubsection*{$c$--strong inequalities}
For $S \subseteq K$ consider the arc capacity inequality written as \eqref{eq:arc-base}.
Relaxing the inequality by dropping $x_i, \ i \in K\setminus S$ and applying integer rounding,
we obtain the so-called \textit{$c$-strong inequality}
\begin{equation}
\label{ineq:c-strong}
\hskip 4mm
\sum_{i \in S} x_i  \leq c_S + y,
\end{equation}
where $c_S =  |S| - \CEIL{a(S) - a_0}$.
The set $S$ is said to be {\it maximal $c$--strong} if
$c_{S \setminus \{i\}} = c_S$ for all $i \in S$ and
$c_{S \cup \{i\}} = c_S +1$ for all $i \in K \setminus S$.
Inequality \iref{ineq:c-strong}
is facet--defining for $conv(\FB)$ if and only if
$S$ is maximal $c$--strong.

Given a point $(\bar x, \bar y)$, there is a $c$--strong inequality violated by $(\bar x, \bar y)$
if and only if there exists a set $ S \subseteq K$ such that
$\sum_{i \in S} \bar x_i - c_{S} > \bar y$.
Then, a $c$--strong inequality is violated if and only if
$\max_{S \subseteq K} \big \{ \sum_{i \in S} \bar x_i -  \FLOOR{a_0 + \sum_{i \in S} (1-a_i)} \big \} =
\max \big \{ \sum_{i \in K} \bar x_i z_i - w: 
\sum_{i \in K} (1-a_i) z_i + a_0 + 1/\lambda \leq w , \ z \in \{0,1\}^K, \ w \in \ZZ \big \} + 1> \bar y$,
where $\lambda$ is the least common multiple of the denominators of the rational numbers $(1-a_i)$ and $a_0$.
The last maximization problem with the constant term $-a_0-1/\lambda$ is \NP--hard. 
Nevertheless, the separation problem
has an optimal solution $(z^*, w^*)$ such that $z^*_i = 1$ if
$\bar x_i = 1$ and $z^*_i = 0$ if $\bar x_i = 0$.
Therefore, we can fix such variables to their optimal values
and solve the separation problem over $i \in K$ such that $0 < \bar x_i < 1$,
which in practice can be done very efficiently even by enumeration, 
as most variables take on values either 0 or 1
in the LP relaxations of network design problems.

\subsubsection*{$k$--split $c$--strong inequalities}
\label{sec:kc}

The $c$-strong inequalities can be generalized by considering a relaxation of the capacity constraints, 
where the integer capacity variables are allowed to take values that are integer multiples of $1/k$
for a positive integer $k$.
Thus the \textit{k-split} relaxation takes the form
\[
\FU^k = \Big \{(x,y) \in \{0,1\}^K \times \ZZ :\sum_{i \in K} a_i x_i \leq a_0 + z/k \  \Big \} \cdot
\]
Letting $c_S^k = \sum_{i \in S} \CEIL{k a_i} - \CEIL{ka(S)-ka_0}$, 
we define the {\it k--split $c$--strong inequality} as
\begin{equation}
\label{ineq:kc}
\sum_{i \in S} \CEIL{ka_i} x_i + \sum_{i \in K \setminus S} \FLOOR{ka_i} x_i \leq c_S^k + ky.
\end{equation}

The $k$--split $c$--strong inequality \eqref{ineq:kc} 
is facet--defining for $conv(\FU)$ if 
$(i)$ $S$ is maximal $c$--strong in the $k$--split relaxation,
$(ii)$ $f_S > (k-1)/k$ and $a_0 \geq 0$, 
$(iii)$ $a_i > f_S$ for all $i \in S$, 
$a_i < 1- f_S$ for all $i \in K \setminus S$,  where
$f_S = a(S)- a_0 - \FLOOR{a(S)-a_0}$.

\begin{example}
	\label{sec:ex-arc}
	Consider the unsplittable arc set
	$$\FU = \bigg \{(x,y) \in \{0,1\}^5 \times \ZZ: \frac{1}{3}x_1+\frac{1}{3}x_2+\frac{1}{3}x_3+\frac{1}{2}x_4+\frac{2}{3}x_5 \leq y \bigg \} \cdot 
	$$
	The maximal $c$--strong inequalities for \FU \ are:
	\[	\begin{array}{ll}
	c_S = 0: &
	x_1 \leq y, \ \
	x_2 \leq y, \ \
	x_3 \leq y, \ \
	x_4 \leq y \\
	c_S = 1: &
	x_1+x_2+x_4 \leq 1 +y, \ \
	x_1+x_2+x_5 \leq 1 +y, \ \
	x_2+x_3+x_4 \leq 1 +y, \ \\
	& x_2+x_3+x_5 \leq 1 +y, \ \
	x_1+x_3+x_4 \leq 1 +y, \ \
	x_1+x_3+x_5 \leq 1 +y \ \\
	c_S = 2: &
	x_1+x_2+x_3+x_4+x_5 \leq 2 +y \\
	\end{array}
	\]
	As the inequalities are maximal, they are facet-defining  for $conv(\FU)$.
	The 2--split $c$-strong inequality
	$x_2+x_3+x_4+x_5 \leq 2y$
	and
	the 3--split $c$-strong inequality 
	$x_1+x_2+x_3+2x_4+2x_5 \leq 3y$ 
	are also facet--defining for $conv(\FU)$.
\end{example}

\subsubsection*{Lifted knapsack cover inequalities}
\label{sec:cover}

Facets different from $c$-strong and $k$-split $c$-strong inequalities can be obtained 
by lifting cover inequalities 
from knapsack restrictions of \FB.
Let $K_0$ and $K_1$ be two disjoint subsets of $K$ and $\FIXY$ be a nonnegative integer. 
Consider the 0-1 knapsack set \FUP \ obtained by restricting  the capacity variable $y$ to \FIXY,
all binary variables indexed with $K_0$ to 0 and
all binary variables indexed with $K_1$ to 1, i.e.,
$\FUP \equiv \{ (x,y) \in \FB: y = \FIXY, \  
x_i = 0 \ \text{for all} \ i \in K_0 \ \text{and} \ x_i = 1 \ \text{for all} \ i \in K_1 \}$.
For this knapsack restriction
$C = K \setminus (K_0 \cup K_1)$ is called a {\it cover} if
$r = a(C) + a(K_1) - a_0 - \FIXY > 0$. 
$C$ is said to be a {\it minimal cover} if $a_i \geq r$ for all $i \in C$.

The cover inequality $\sum_{i \in C} x_i \leq |C|-1$ is
facet--defining for $conv(\FUP)$ if and only if $C$ is a minimal cover. 
One practical way of lifting inequalities
is sequential lifting, in which restricted variables are introduced to an inequality one at a time
in some sequence. 
A lifted cover inequality 
\begin{equation}
\label{ineq:lifted}
\sum_{i \in C} x_i + \sum_{i \in K_0} \alpha_i x_i + \sum_{i \in K_1} \alpha_i (1-x_i) + \alpha (\FIXY-y) \leq |C|-1
\end{equation}
can be constructed in $O(|K|^3)$ if
the capacity variable $y$ is lifted first and  such inequalities subsume all
$c$--strong inequalities.

\begin{example}
	\label{sec:ex}
	For $\FU$ given in Example~\ref{sec:ex-arc}
	we list below the lifted cover inequalities of \FU \ that are not $c$--strong inequalities.
	\vskip -1mm
	\begin{align*}
	\hskip -2mm
	\begin{array}{r | c | l}
	\FIXY  \hskip 1mm & (C, \ K_0, \ K_1) & \hskip 2cm \text{Inequalities} \\
	\hline
	1  \hskip 1mm &  \hskip 1mm (\{2,3,4\}, \{1,5\}, \emptyset)  \hskip 1mm &  
	\hskip 1mm x_2+x_3+x_4+x_5 \leq 2y \\
	1  \hskip 1mm &  \hskip 1mm (\{1,4,5\}, \{2,3\}, \emptyset)  \hskip 1mm &  
	\hskip 1mm x_1+x_2+x_4+x_5 \leq 2y \ \text{and} \
	x_1+x_3+x_4+x_5 \leq 2y \\
	2  \hskip 1mm &  \hskip 1mm (\{1,2,3,4\}, \emptyset, \{5\})  \hskip 1mm &  
	\hskip 1mm x_1+x_2+x_3+x_4+2x_5 \leq 2y+1 \\
	2  \hskip 1mm &  \hskip 1mm (\{1,2,3,5\}, \emptyset, \{4\})  \hskip 1mm &  
	\hskip 1mm x_1+x_2+x_3+2x_4+x_5 \leq 2y+1
	\end{array}
	\end{align*}
\end{example}

\ignore{
	Computational results in \cite{BGW:private} suggest that $c$-strong inequalities are quite effective in solving unsplittable multi-commodity network design problems.
	Computational studies in \cite{AR:arc} and \cite{HKLS:edge-cap} indicate that while the $k$-split $c$-strong and the lifted knapsack cover inequalities provide additional strengthening of the relaxations, the marginal impact on top of the basic $c$-strong inequalities is limited. The latter result may be due to the lack of efficient separation procedures for these inequalities.
}
Computational results suggest that $c$-strong inequalities are quite effective in solving unsplittable multi-commodity network design problems.
Moreover, while the $k$-split $c$-strong and the lifted knapsack cover inequalities provide additional strengthening of the relaxations, the marginal impact on top of the basic $c$-strong inequalities is limited. The latter result may be due to the lack of efficient separation procedures for these inequalities.

\subsection{Multi-facility Arc Set}

\label{sec:arcset-mf}

In this section we consider the multi-facility extension of the arc sets discussed in the previous sections.
\cite{DGW:cont-knap} study a 
continuous knapsack set with two integer capacity variables:
\[
F_{S2} = \bigg \{(x,y) \in [0,1]^K \times \Z_+^2: \sum_{i \in K} a_i x_i \le a_0 + c_1 y_1 + c_2 y_2 \bigg \} \cdot
\]
They show that all non-trivial facet-defining inequalities of $\conv(F_{S2})$ are of the form
\[
\sum_{i \in S} a_i x_i + \gamma_1 y_1 + \gamma_2 y_2 \ge \beta,
\]
where $S \subseteq K$ and $w + \gamma_1 y_1 + \gamma_2 y_2 \ge \beta$ is facet-defining for the set
\[
Q(a,b) = \conv \big \{ (w,y) \in \RR \times \ZZ_+^2: w + c_1 y_1 + c_2 y_2 \ge b, \ a \ge w  \big \}.
\]
It turns out that all facets of $Q(a,b)$ can be enumerated in polynomial time. 
Therefore, for each $S \subseteq K$ non-trivial facet-defining inequalities for $F_{SM}$ can be obtained from relaxations of the form $Q(a,b)$.

For the general case with many facility types
\[
F_{SM} = \bigg \{(x,y) \in [0,1]^K \times \Z_+^M: \sum_{i \in K} a_i x_i \le a_0 + \sum_{m \in M} c_m y_m \bigg \}
\]
a similar approach of complementing the flow variables for a subset $S \subseteq K$, scaling the base inequality by $c_s, \ s \in M$, and applying mixed-integer rounding gives for $F_{SM}$ the \textit{multi-facility residual capacity inequalities}
\[
\sum_{m \in M} \phi_s(c_m) y_m + \sum_{i \in S} a_i(1-x_i) \ge r_s \eta_s, 
\] 
where $r_s = a(S) - a_0 - \floor{(a(S) - a_0)/c_s}c_s$, $\eta_s = \ceil{(a(S) - a_0)/c_s}$,  and for $k \in \ZZ$
\[
\phi_{s}(c) =
\begin{cases}
c - k(c_s-r_{s}) & \text{if} \ kc_s \leq c < kc_s + r_{s}, \\
(k+1)r_{s}       & \text{if} \ kc_s + r_{s} \leq c < (k+1)c_s.
\end{cases}
\]

\section{Valid Inequalities From Cut Sets}
\label{sec:cutset}

In this section we review valid inequalities for the network design problem based on relaxations formed over cuts of the network. 
We first start with the single-facility case and then generalize the inequalities 
for multiple facilities.

\subsection{Single-facility Case}

\begin{figure}
	\centering
	\includegraphics[scale=0.5]{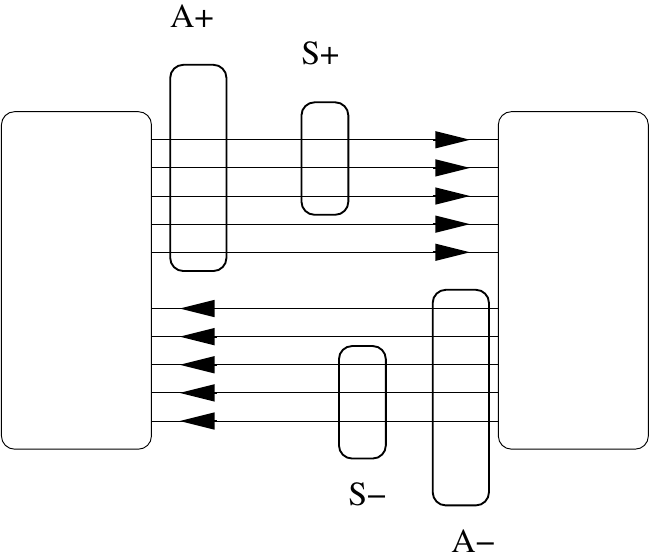}
	\caption{Cut-set relaxation of the network design problem.}
	\label{fig:cut-arcs}
\end{figure}

Consider a nonempty two-partition $(U,V)$ of the vertices of the network.
Let $b^k$ denote the total supply of commodity $k$ in $U$ for $V$.
Let $A^+$ be the set of arcs directed from $U$ to $V$,
$A^-$ be the set of arcs directed from $V$ to $U$, and $A = A^+ \cup A^-$, as shown in Figure~\ref{fig:cut-arcs}.
As before, 	 $x^k_a$ denotes the flow of commodity $k$ on	arc $a\in A$ for $k \in K$. 	
The constraints of the 
multicommodity network design problem across the
cut are
\begin{alignat}{2}
\label{eq:nd:sbalance}
x^k(A^+) - x^k(A^-) & =  b^k, &  \ \ &  k \in K,\\
\label{eq:nd:svub}
\sum_{k \in K} x^k_a & \leq  \bar c_a + c y_a,  & \ \ &  a \in A.
\end{alignat}
Then the corresponding multicommodity cut--set polyhedron is defined as
\begin{alignat}{2}
\nonumber
\FMS = conv \Big \{(x,y) \in  \RR^{A \times K}_+ \times \ZZ^{A}_+ : (x,y) \ \mbox{satisfies \eqref{eq:nd:sbalance} and \eqref{eq:nd:svub}}  \Big \} \cdot
\end{alignat}
We refer to the single-commodity case of \FMS as \FSS.

In the following sections we describe valid inequalities for \FMS 
by considering single-commodity relaxations of $\FMS$ obtained by
aggregating flow variables and balance equations \eqref{eq:nd:sbalance} over subsets of $K$.
For $Q \subseteq K$ let $x_Q(S) = \sum_{k \in Q} x^k(S)$ and $b_Q = \sum_{k \in Q} b^k$.

\vskip 3mm \noindent
\textit{Cut-set inequalities.}
Consider the following relaxation of \FMS on the integer capacity variables:
\begin{align*}
\update{\bar{c}(A^+) +} c y(A^+) \ge	x_K(A^+)  & \ge   b_K.
\end{align*}
Applying the integer rounding procedure reviewed in Section~\ref{sec:mir} to this relaxation, one obtains the so-called \textit{cut-set inequality}
\begin{align}
\label{eq:nd:cut-set}
y(A^+) \geq \CEIL{(b_K \update{- \bar{c}(A^+)})/c} 
\end{align}
for \FMS, which is unique per cut-set relaxation. 
Finding the best cut-set relaxation is not easy however. For
the single-source single-sink case, the problem of finding the best cut set can be posed
as an $s-t$ max-cut problem. 

\vskip 3mm \noindent
\textit{Flow-cut-set inequalities.}
The basic cut-set inequalities \eqref{eq:nd:cut-set} can be generalized 
by incorporating the flow variables in addition to the capacity variables (see Figure~\ref{fig:cut-arcs}).	
For $S^+ \subseteq A^+$, $S^- \subseteq A^-$ and $Q \subseteq K$ consider the following relaxation of \FMS:
\begin{align*}
\update{\bar c(S^+) + }	c y(S^+) +  x_Q(A^+ \setminus S^+)  - x_Q(S^-) & \ge   b_Q, \\
0 \leq \sum_{k \in Q} x^k_a & \leq  \update{\bar c_a + }  c y_a,  \ \forall a \in A, 
\end{align*}
which is written equivalently as
\begin{align*}
c \big [y(S^+)- y(S^-) \big] +  x_Q(A^+ \setminus S^+) + \big [ \bar c(S^-) + c y(S^-) - x_Q(S^-) \big] & \ge   
\update{b_Q'}, \\
0 \leq \sum_{k \in Q} x^k_a & \leq  \update{\bar c_a + } c y_a,  \ \forall a \in A. 
\end{align*}
where \update{$b'_Q =  b_Q -\bar c(S^+) + \bar c(S^-)$}.
Letting $r_Q = \update{b_Q'} - \FLOOR{\update{b_Q'}/c}c$ and $\eta_Q = \CEIL{\update{b_Q'}/c}$ and 
observing that  $x_Q(A^+ \setminus S^+)  \ge 0$, $\bar c(S^-) + c y(S^-) - x_Q(S^-) \ge 0$,
we can apply the MIR procedure reviewed in Section~\ref{sec:mir} to this relaxation to arrive at 
the flow-cut set inequalities 
\begin{align}
\label{eq:nd:flow-cut-set-inout}
r_Q y(S^+)+ x_Q(A^+ \setminus S^+) + (c-r_Q) y(S^-) - x_Q(S^-) \geq r_Q \eta_Q .
\end{align}

The flow-cut-set inequalities \eqref{eq:nd:flow-cut-set-inout} along with the balance, bound, and capacity constraints are sufficient to describe the single-commodity case \FSS.

For a given $Q \subseteq K$, observe that flow-cut-set inequalities~\eqref{eq:nd:flow-cut-set-inout} is an exponential class. 
However, given a solution $(\bar x, \bar y)$,
one finds a subset $S^+$ with the smallest left-hand-side value as follows:
if $r_Q \bar y_a < \sum_{k \in Q} \bar x^k_a$ for $a \in A^+$, then
we include $a$ in $S^+$; if $(c-r) \bar y_a < \sum_{k \in Q} \bar x^k_a$ for $a \in A^-$, then we include $a \in S^-$.

For a fixed cut of the network,
the complexity of separating multi-commodity flow cut--set inequalities \eqref{eq:nd:flow-cut-set-inout}
is an open question.
Optimization of a linear function over \FMS \ is \NP--hard as the facility location
problem is a special case of it.
For a multi-commodity single-facility network design problem of a single arc,
cut--set inequalities \eqref{eq:nd:lm-cutset} reduce to the
residual capacity inequalities \eqref{ineq:rc}, for
which an exact linear--time separation method is known. 
From here it follows that, for a single-facility problem, if $S^+$ and $S^-$ are fixed,
then one can find a subset of commodities $Q \subseteq K$ that gives a most violated
inequality in linear time. Alternatively, if $Q$ is fixed, since the model
reduces to a single-commodity, one can find the subsets $S^+ \subseteq A^+$ and 
$S^- \subseteq A^-$ that give a most violated inequality in linear time as well.
However, the complexity of determining $Q, \ S^+$ and $S^-$ simultaneously is an open question.

\begin{example}
	Consider the following single-commodity optimization problem with two outflow arcs and one inflow arc:	
	\[
	\max \ x_1+x_2+x_3-y_1-y_2-y_3 \ \st \  x_1+x_2-x_3=0.5, \ 0 \le x_i \le y_i \in \ZZ, \ i=1,2,3.
	\]
	One of its fractional solutions is
	$x_1=y_1=0.5$ and all other variables zero. Adding the cut-set inequality 
	$$y_1 + y_2 \ge 1$$ 
	cuts off this solution, but leads to another fractional solution:
	$x_1=y_1=1, \ x_3=y_3=0.5$
	Adding the flow-cut-set inequality 
	$$0.5x_1+y_2+0.5x_3-y_3 \ge 0.5$$
	cuts it off, but gives the fractional solution: 
	$x_2=y_2=1, x_3=y_3=0.5$.
	Adding the flow-cut-set inequality 
	$$y_1+0.5x_2+0.5x_3-y_3 \ge 0.5$$
	cuts it off, but this time gives the fractional solution: 
	$x_1 = y_1 = x_2=y_2= x_3=y_3=0.5$.
	Finally, adding the flow-cut-set inequality 
	$$0.5x_1+0.5x_2+0.5x_3-y_3 \ge 0.5$$
	leads to an optimal integer solution $x_1=0.5, y_1=1$ and all other variables zero. 
\end{example}

\subsection{Multi-facility Case}
\label{sec:mm}

Next we consider network design models
where one is allowed to install facilities of multiple types
with different capacities on the arcs of the network. 
Let $c_m$ be the capacity of facility of type $m, \ m \in M$.
No assumption is made on either the number of facility types or 
the structure of capacities (other than $c_m > 0$ and rational).

\ignore{
	\cite{MM:st-path,MMV:modeling-2fnlp,PW:div,BG:cnd,G:bc-nd,WY:contknap-div} give valid inequalities
	for the the network design problem with multiple capacities when capacities are divisible.
	\cite{ANS:avub} consider a binary capacity version with no assumption on divisibility.
	Multi-commodity multi-facility network
	design problems are considered in \cite{BCGT:cap-inst,A:nd}.	
}

The constraints of the 
multi-commodity  multi-facility  problem across cut $A$ are written as
\begin{alignat}{2}
\label{eq:nd:mbalance}
x^k(A^+) - x^k(A^-) & =  b^k, &  \ \ & ~~~ k \in K,\\
\label{eq:nd:mavub}
\sum_{k \in K} x^k_a & \leq \update{\bar c_a + } \sum_{m \in M} c_m y_{m,a},  & \ \ & ~~~ a \in A.
\end{alignat}
So the corresponding multi-commodity multi-facility cut-set polyhedron is 
\begin{alignat}{2}
\nonumber
\FMM = conv \Big \{(x,y) \in  \RR^{A \times K}_+ \times \ZZ^{A \times M}_+ : (x,y) \ \mbox{satisfies \eqref{eq:nd:mbalance} and \eqref{eq:nd:mavub}} \Big \}.
\end{alignat}

For $Q \subseteq K$ and $s \in M$ let $r_{s,Q} = \update{b_Q'} - \FLOOR{\update{b_Q'} / c_s}c_s$,
$\eta_{s,Q} = \CEIL{\update{b_Q'} / c_s}$. Then, the following {\it multi-commodity multi-facility cut-set inequality} is valid for \FMM:
\begin{equation}
\label{eq:nd:lm-cutset}
\sum_{m \in M} \phi_{s,Q}^+(c_m) y_m(S^+) \hspace{-0.5mm} + 
\hspace{-0.5mm} x_Q(A^+ \hspace{-0.5mm} \setminus S^+)  \hspace{-0.5mm} + 
\hspace{-1mm} \sum_{m \in M} \phi_{s,Q}^-(c_m) y_m(S^-) \hspace{-0.5mm} - \hspace{-1mm} x_Q(S^-) \hspace{-0.5mm} \geq 
\hspace{-0.5mm} r_{s,Q} \eta_{s,Q},
\end{equation}
where, for $k \in \ZZ$, 
\[
\phi_{s,Q}^+(c) =
\begin{cases}
c - k(c_s-r_{s,Q}) & \text{if} \ kc_s \leq c < kc_s + r_{s,Q}, \\
(k+1)r_{s,Q}       & \text{if} \ kc_s + r_{s,Q} \leq c < (k+1)c_s,
\end{cases}
\]
and
\[
\phi_{s,Q}^-(c) =
\begin{cases}
c - kr_{s,Q}       & \text{if} \ kc_s \leq c < (k+1)c_s - r_{s,Q}, \\
k(c_s-r_{s,Q})     & \text{if} \ kc_s - r_{s,Q} \leq c < kc_s. \\
\end{cases}
\]
The above $\phi_{s,Q}^+$ and $\phi_{s,Q}^-$ are subadditive MIR functions written in closed form. The
multi-commodity  multi-facility cut-set inequality \eqref{eq:nd:lm-cutset} is facet-defining for $\FMM$ if
$(S^+, A^+ \setminus S^+)$ and $(S^-, A^- \setminus S^-)$ 
are nonempty partitions, $r_{s,Q} > 0$, and $b^k > 0$ for all $k \in Q$. 

For the single-commodity case \FSM, inequalities \eqref{eq:nd:lm-cutset} reduce to
\begin{equation}
\label{eq:nd:sm-cutset}
\sum_{m \in M} \phi_s^+(c_m) y_m(S^+) + x(A^+ \setminus S^+) +
\sum_{m \in M} \phi_s^-(c_m) y_m(S^-) - x(S^-) \geq r_s \eta_s.
\end{equation}
In this case, given a cut $A$ 
for each facility $s \in M$, the multi-facility cut-set
inequalities \eqref{eq:nd:sm-cutset} is an exponential class by the choice 
of the subsets of arcs $S^+$ and $S^-$. 
However, finding a subset that gives a most
violated inequality for a point $(\bar x, \bar y)$ is straightforward.
If $\sum_{m \in M} \phi_s^+(c_m) \bar y_{m,a} < \bar x_a$ for $a \in A^+$, 
then we include $a$ in $S^+$, and 
if $\sum_{m \in M} \phi_s^-(c_m) \bar y_{m,a} < \bar x_a$ for $a \in A^-$,
then we include $a$ in $S^-$. 
Since $\phi_s^+(c)$ or $\phi_s^-(c)$ can be calculated in constant time,
for a fixed cut $A$, a
violated multi-facility cut-set inequality is found in $O(|A||M|)$
if there exists any.

\begin{example}
	We specialize inequality \eqref{eq:nd:lm-cutset} for the network
	design problem with two types of facilities. 
	Let the vectors $y_1$ and $y_2$ denote the facilities with capacities
	$c_1 = 1$ and $c_2 = \lambda >1$ with $\lambda \in \ZZ$, respectively. 
	Let $Q$ be a nonempty subset of the commodities. 
	Then by letting $s = 1$, we have $r_{1,Q} = b_Q - \FLOOR{b_Q}$ and 
	inequality \eqref{eq:nd:lm-cutset} becomes
	\begin{multline*}
	r_{1,Q}y_1(S^+)+(r_{1,Q} \FLOOR{\lambda}+\min\{\lambda-\FLOOR{\lambda},r_{1,Q}\})y_2(S^+) + 
	x_Q(A^+ \setminus S^+)+ 
	(1-r_{1,Q}) \\ y_1(S^-)+((1-r_{1,Q})\FLOOR{\lambda}+\min\{\lambda-\FLOOR{\lambda},1-r_{1,Q}\})y_2(S^-) -
	x_Q(S^-) \geq r_{1,Q} \CEIL{b_Q},
	\end{multline*}
	which, when $\lambda$ is integer, reduces to 
	\begin{multline}
	\label{ex2}
	r_{1,Q}y_1(S^+)+\lambda r_{1,Q}y_2(S^+) + x_Q(A^+ \setminus S^+)+ \\
	(1-r_{1,Q})y_1(S^-)+\lambda(1-r_{1,Q})y_2(S^-)-x_Q(S^-) \geq r_{1,Q} \CEIL{b_Q}.
	\end{multline}
	Notice that inequality \eqref{ex2} is not valid for $\FMM$ unless $\lambda \in \ZZ$.
	Also by letting $s=2$, we have $r_{2,Q}=b_Q-\FLOOR{b_Q/\lambda}\lambda$. 
	So the corresponding multi-commodity two-facility cut-set inequality is
	\begin{multline*}
	\label{ex3}
	\min\{1,r_{2,Q}\}y_1(S^+)+r_{2,Q}y_2(S^+)+x_Q(A^+ \setminus S^+)+ \\
	\min\{1,\lambda-r_{2,Q}\}y_1(S^-)+
	(\lambda-r_{2,Q})y_2(S^-)-x_Q(S^-) \geq r_{2,Q} \CEIL{b_Q/\lambda}.
	\end{multline*}		
\end{example}

It should be clear that multi-facility flow-cut-set inequalities are also valid for a single-facility model with varying capacities on the arcs of the network.  \ignore{Moreover, the existing arc capacities $\bar c_a$, $a \in A$, can be treated through a relaxation by introducing auxiliary capacity variables $\bar y_a$ with capacity unit $\bar c_a$, which are fixed to one after writing an inequality. {\cred Due to sub-additivity, isn't this a weakening of the actual MIR inequality? Why not define $b_Q$ as demand minus existing $A^+$ capacity plus existing $A^-$ capacity?}}
The inequalities from cut-set relaxations have been shown to be very effective in solving network design problems in computational studies. 


\section{Partition Inequalities}
\label{sec:part}

Partition inequalities are arguably the most effective cutting planes for the network design problem. These inequalities have non-zero coefficients for only the integer capacity variables 
that cross a multi-cut obtained from a partition of the nodes of the network. 
They generalize the cut-set inequality \iref{eq:nd:cut-set} described in Section~\ref{sec:cutset}.

For ease of exposition, first  consider a two-partition of the node set $\cN=N_1\cup N_2$.
As discussed in Section \ref{sec:prelim-part} shrinking node sets $N_1$ and $N_2$ leads to a network with two nodes and two edges (assuming there is at least one arc from a node in $N_1$ to a node in $N_2$, and vice versa). Then, the inequality 
\beq \sum_{m\in M}c_my_{m12}\ge \sum_{k \in K}\sum_{v \in N_2}w^k_v-\sum_{u\in N_1}\sum_{v \in N_2} \bar c_{uv}~~=~b\eeq
must be satisfied by all feasible solutions of the two-node problem. Following the argument in
Section \ref{sec:prelim-part}, inequality
\beqn \sum_{m\in M}c_m \bigg (\sum_{u\in N_1}\sum_{v \in N_2}y_{muv} \bigg )\ge b\label{ineq.cutset.base}\eeqn
is valid for the solutions to the (LP relaxation of the) original problem.
Notice that  inequality \iref{ineq.cutset.base} is a metric inequality \iref{eq:metax} generated  by the vector $ v\in \{0,1\}^{A}$, where $v_{ij}=1$ if and only if $i\in N_1$ and $j\in N_2$.

As we assume that all $c_m$ are integral, which is the case in most applications, the inequality \iref{ineq.cutset.base} leads to the integer knapsack cover set 
$$X=\bigg\{z\in\Z^M:\sum_{m\in M}c_mz_m\ge b,~ z\ge 0\bigg\},$$  
where the variable $z_m$ stands for the sum $\sum_{u\in N_1}\sum_{v \in N_2}y_{muv}$.
Consequently, any valid inequality $\sum_{m\in M}\alpha_mz_m \ge \beta$ for $X$ yields a valid inequality for $\poly$ after replacing each variable $z_m$ with the corresponding sum of the original variables.

The polyhedral structure of the set $X$ when $c_{m+1}$ is an integer multiple of $c_m$ for all $m=1,\ldots,|M|-1$ has been studied by \cite{PW:div} who derive what they call ``partition" inequalities and show that these inequalities  together with the nonnegativity constraints describe $\conv(X)$.
They also derive conditions under which partition inequalities are valid in the general case when the divisibility condition does not hold.

The partition inequalities described in \cite{PW:div} 
are obtained by applying the MIR procedure iteratively.
More precisely, the first step is to chose a subset $\{j_1,j_2,\ldots,j_r\}$ of the index set $M$, where $j_i<j_{i+1}$, and therefore, $c_{j_i}<c_{j_{i+1}}$ for all $i=1,\ldots,r-1$.
The inequality $\sum_{m\in M}c_mz_m\ge b$ is then divided by $c_{j_r}$ and the MIR cut based on this inequality is written. 
The resulting MIR inequality is then divided by $c_{j_{r-1}}$ and the MIR procedure is applied again.
This process is repeated with all $c_{j_i}$ for  $i=1,\ldots,r$ to obtain the final inequality.
Note that the sequential application of the MIR procedure yields valid inequalities even when the divisibility condition does not hold. 
However, in this case, they are not sufficient to define $\conv(X)$.


Now consider a three-partition of the node set $\cN=N_1\cup N_2\cup N_3$.
Following the discussion on two-partitions, consider the single capacity network design problem with 
$ \tilde\cG =( \tilde\cN, \tilde\cA)$ where $ \tilde\cN=\{1,2,3\}$ and $\tilde\cA=\{a_{12},a_{13},a_{21},a_{23},a_{31},a_{32}\}$.
Furthermore, let $\bar c_a$ and $\bar t_a$ respectively denote the existing capacity and traffic demands for $a\in\tilde\cA$.
Clearly, any valid inequality for the simplified problem on $\tilde G$ can be transformed into a valid inequality for the original problem defined on $ \cG$.

For the three-node problem, there are three possible two-partitions and for each partition, one can write two possible cut-set inequalities by treating one of the two sets as $N_1$ and the other as $N_2$. 
Consequently, one can write six different two-partition inequalities where each capacity variable appears in exactly two inequalities. Summing all six inequalities leads to 
\beqn \sum_{a\in \tilde A} \sum_{m\in M^{\emptybox{.13cm}}}2 c_m y_{ma}~\ge~\sum_{ i\in \tilde \cN}\ceil{s_i}+\sum_{ i\in \tilde \cN}\ceil{t_i}, \label{ineq:sum2part}
\eeqn
where $s_i$ denotes the difference between the traffic leaving node $i$ and the existing capacity on the outgoing arcs.
Similarly $t_i$ is the difference between the traffic entering node $i$ and the existing capacity on the incoming arcs.
If the right-hand-side of inequality \eqref{ineq:sum2part} is an odd number, dividing the inequality by two and rounding up the right-hand-side yields the following inequality
\beqn \sum_{m\in M^{\emptybox{.13cm}}}c_m\sum_{a\in \tilde A}  y_{ma}~\ge~\xceil{	\frac{\sum_{ i\in \tilde \cN}\ceil{s_i}+\sum_{ i\in \tilde \cN}\ceil{t_i})}2}.\label{ineq:3part}\eeqn

Next we will generate a similar inequality based on metric inequalities.
Let $a,b\in\tilde{\cN}$ be two distinct nodes and define $v^{ab}\in \{0,1\}^{\tilde\cA}$ to be the vector with components $v_{ab}=v_{ac}=v_{bc}=1$ and $v_{ba}=v_{ca}=v_{cb}=0$. 
Now consider the metric inequality \iref{eq:metax} generated  by $v^{ab}$:
\beq \sum_{m\in M}c_m (y_{mab}+y_{mac}+y_{mbc})~\ge~\bar t_{ab}+\bar t_{ac}+\bar t_{bb}-\bar c_{ab}-\bar c_{ac}-\bar c_{bc}.\eeq
Once again, if fractional, the right-hand-side of this inequality can be rounded up. 
In addition, more valid inequalities can be generated using the MIR procedure iteratively.

Furthermore, let $c\in\tilde\cN$ be the node different from $a$ and $b$ and note that adding up the metric inequalities  generated  by $v^{ab}$ and $v^{cb}$, one obtains
\beqn \sum_{m\in M^{\emptybox{.153cm}}}c_m \sum_{a\in \tilde A} y_{ma}~\ge~\ceil{d_{ab}}+\ceil{d_{cb}}\label{ineq:3part.metric}\eeqn
where $d_{ij}$ denotes the right-hand-side of the metric inequality generated by  $v^{ij}$ for $(i,j)\in\cA$.
Moreover, adding up the metric inequalities  generated  by $v^{ba}$ and $v^{ca}$ gives a similar inequality with right-hand-side of $\ceil{d_{ba}}+\ceil{d_{ca}}$. 
Similarly, $v^{ac}$ and $v^{bc}$ yields an inequality with right-hand-side of $\ceil{d_{ac}}+\ceil{d_{bc}}$.
Adding up two of these inequalities with the larger right-hand-side and dividing the resulting inequality by two leads to a valid inequality of the form \iref{ineq:3part}. 
More precisely, if both $\ceil{d_{ba}}+\ceil{d_{ca}}$ and  $\ceil{d_{ab}}+\ceil{d_{cb}}$ are larger than $\ceil{d_{ac}}+\ceil{d_{bc}}$, then the resulting inequality is
\beqn \sum_{m\in M^{\emptybox{.153cm}}}c_m \sum_{a\in \tilde A} y_{ma}~\ge~\xceil{	\frac{\ceil{d_{ba}}+\ceil{d_{ca}}+\ceil{d_{ab}}+\ceil{d_{cb}} }{2}}\label{ineq:3part.metric.CG}.\eeqn

%
%

In addition to  inequalities \iref{ineq:3part} and \iref{ineq:3part.metric.CG}, it is possible to write similar {\em total capacity inequalities} 
by combining some cut-set inequalities with metric inequalities in such a way that the left-hand-side of the inequality has all the capacity variables with a coefficient of two.
As all these inequalities have the same left-hand-side, only the one with the largest right-hand-side should be used.
For example, if $\bar t_{ij}=1/2$ and $\bar c_{ij}=0$ for $(i,j)\in\tilde \cA$, then the right-hand-side of  \iref{ineq:3part} is 3, whereas the right-hand-side of  \iref{ineq:3part.metric.CG} is 4 and therefore inequality  \iref{ineq:3part.metric.CG} is stronger than \iref{ineq:3part}.
However, if $\bar t_{ij}=1/3$ and $\bar c_{ij}=0$ for $(i,j)\in\tilde \cA$, then the right-hand-side of  \iref{ineq:3part} is still 3, whereas the right-hand-side of  \iref{ineq:3part.metric.CG} becomes 2.

Furthermore, as  total capacity inequalities  have the same form as inequality \iref{ineq.cutset.base}, one can define the corresponding integer knapsack cover set from the stronger one and derive further valid inequalities using the MIR procedure iteratively.

\cite{A:Hamid:2015} study the undirected variant of the two-facility network design problem, where the total flow on an arc plus the flow on reverse arc is limited by the capacity of the undirected edge associated with the two arcs. 
In this case, the authors computationally enumerate the complete list of facets that can be obtained from a given three-partition.
Also see \cite{A:kpart} for a study of four-partition facets for the undirected variant of the single-facility network design problem.
Their computational study suggests that using larger partitions of the node set improves the relaxation but with diminishing returns. 
More precisely,  they observe that two, three and four-partition cuts reduce the optimality gap of the LP relaxation to 12.5\%, 6.3\%, and 2.6\%, respectively.



\section{Bibliographical Notes}
\subsubsection*{Introduction  }
Capacity expansion problems have been studied in the context of telecommunucation \citep{MW:nd,M:netsyn,BMSW:nexp,BMW:local-access} and train scheduling with multiple locomotive options \citep{FBFGN:train}.
Unsplittable flow version of the network design problem appears in telecommunication networks using multiprotocol label switching (MPLS) technique,
production and distribution with single sourcing, and express package delivery, train scheduling \citep[e.g.][]{GA:backbone,BHV:intmulticom,JP:train}.  
\cite{OKY:splitdelivery} utilize the splittable flow model to solve the split delivery vehicle routing problem.
\citet{AG:capacity} study the connection between valid inequalities for directed, bidirected and undirected network design problems.
	
\subsubsection*{Problem formulation and preliminaries  }
\cite{DBLP:journals/siamdm/Gunluk07} describes how to obtain the smallest number of commodities necessary to formulate a given multicommodity flow problem. 
Metric inequalities and their extensions have been used for various network design problems by several authors, including 
\cite{DahlStoer1998,M:proj-netload,LY:project,AMS:metric,CCG:mc-netdes,BCGT:cap-inst}. 
\cite{M:sep-metric} presents computations that illustrate the value of utilizing metric inequalities through a bi-level programming separation procedure. \cite{A:Hamid:2015} show that if inequality \iref{ineq_small} is facet-defining for the network design problem on $\tilde\cG$ with $\tilde c$ and $\tilde w$, then inequality \iref{ineq_big} is facet-defining for \poly\ provided that $\tilde\alpha=0$,~ $\gamma> 0$, and node sets $N_1,\ldots,N_p$ induce connected components of  $G$.
In addition, \cite{RKOW:cutset} show that the same  result holds without the assumption $\tilde\alpha=0$  when $p=2$ and $|M|=1$.
\subsubsection*{Valid inequalities from arc sets}
The arc sets and their generalizations are studied by 
\cite{MMV:conv-2core,AR:arc,HKLS:edge-cap,BGW:private,AG:mixset,Y:arcset}.
\cite{MMV:conv-2core} introduced the {\it residual capacity inequality}. 
\cite{AR:arc} gave the polynomial-time separation procedure for residual capacity inequities.
\citet{BGW:private} introduced the \textit{$c$-strong inequality}
\eqref{ineq:c-strong} for \FB. The
{$k$-split $c$-strong inequalities}
\eqref{sec:kc} are given by \cite{AR:arc}.
\citet{AC:2int} show that all facets of $Q(a,b)$ can be enumerated in polynomial time.
\subsubsection*{Valid inequalities from cut sets}
A recent review on cut-based inequalities for network design can be found in
\citet{RKOW:cutset}. \citet{MM:st-path} introduce the integer cut-set inequalities for \FSS. 
For the single-source single-sink case, given a solution,  \citet{B:nd-cut} formulated the problem of finding the best cut set as an $s-t$ max-cut problem. 
\citet{BG:cnd,CGS:nd-batch} gave the mixed-integer flow-cut-set generalization \eqref{eq:nd:flow-cut-set-inout} of the basic cut-set inequalities.
\citet{A:nd} showed that the flow-cut-set inequalities \eqref{eq:nd:flow-cut-set-inout} along with the balance, bound, and capacity constraints are sufficient to describe the single-commodity case \FSS.
	\cite{MM:st-path,MMV:modeling-2fnlp,PW:div,BG:cnd,G:bc-nd,WY:contknap-div} give valid inequalities
for the  network design problem with multiple capacities when capacities are divisible.
\cite{ANS:avub} consider a binary capacity version with no assumption on divisibility.
Multi-commodity multi-facility network
design problems are considered in \cite{BCGT:cap-inst,A:nd}.

\subsubsection*{Partition inequalities}
The polyhedral structure of the set $X$  without the divisibility assumption on the capacities has been studied by \cite{A:mip} and \cite{Y:knap-cover}, also see the references therein.
\citet{AGK:3part-fc} gave three-partition flow cover inequalities for the fixed-charge network design problem. The total capacity inequalities were proposed by \cite{BCGT:cap-inst}. 
\cite{A:Hamid:2015} has enumerated the complete list of facets that can be obtained from a given three-partition for the undirected variant of the two-facility network design problem.
\cite{A:kpart} has studied the four-partition facets and computationally showed that using larger partitions of the node set improves the LP relaxation but with diminishing returns.

\section{Conclusions and Perspectives}
\label{sec:conc}

In this chapter, we reviewed strong  valid inequalities for the multi-commodity, multi-facility network design problem. Metric inequalities for projecting out continuous flow variables, mixed-integer rounding from appropriate base relaxations, and shrinking the network to a small $k$-node graph have been the main tools for deriving the inequalities introduced in the literature. Going forward, we expect more recent techniques such as multi-step mixed-integer rounding \citep{DG:twostep}, mingling \citep{AG:mingle}, and multi-step mingling \citep{AK:nstepmingle} that generalize and extend mixed-integer rounding to be useful for deriving new inequalities for this class of network design problems with continuous as well as general integer variables.

We finish this section with comments on computational aspects and challenges in implementing
branch-and-cut algorithms. One of the prerequisites of applying the valid inequalities described in this chapter as cutting planes is to automatically recognize the appropriate network structure and identify flow-balance constraints and the relevant capacity constraints as part of a hidden multi-commodity, multi-facility network structure. \cite{AR:mcf-sep} describe algorithms to automatically detect such network structures and generate inequalities from cut-set relaxations. They utilize c-MIR inequalities \citep{MW:MIR}, which are obtained by first complementing bounded variables and then applying the MIR procedure. Given a fractional solution, it is nontrivial to decide which aggregations to apply to flow variables, which variables to complement. Nevertheless, the automatic detection and separation procedure in this paper results in 18\% time reduction for a large set of publicly available test problems.

\section*{Acknowledgements}
A. Atamt\"urk is supported, in part,
by grant FA9550-10-1-0168 from the Office of the Assistant Secretary of Defense for Research and Engineering.
O. G\"unl\"uk was hosted by the Numerical Analysis Group at Oxford Mathematical Institute and by the Alan Turing Institute during this project and would like to thank the members of both groups for their hospitality. 

\bibliographystyle{apalike}
\bibliography{master}

\newcommand{\SortNoop}[1]{}
\begin{thebibliography}{}

\bibitem[Achterberg and Raack, 2010]{AR:mcf-sep}
Achterberg, T. and Raack, C. (2010).
\newblock The {MCF}-separator: detecting and exploiting multi-commodity flow
  structures in {MIPs}.
\newblock {\em Mathematical Programming Computation}, 2:125--165.

\bibitem[Agarwal, 2006]{A:kpart}
Agarwal, Y.~K. (2006).
\newblock K-partition-based facets of the network design problem.
\newblock {\em Networks}, 47:123--139.

\bibitem[Agra and Constantino, 2006]{AC:2int}
Agra, A. and Constantino, M. (2006).
\newblock Description of 2-integer continuous knapsack polyhedra.
\newblock {\em Discrete Optimization}, 3:95 -- 110.

\bibitem[Atamt{\"u}rk, 2002]{A:nd}
Atamt{\"u}rk, A. (2002).
\newblock On capacitated network design cut--set polyhedra.
\newblock {\em Mathematical Programming}, 92:425--437.

\bibitem[Atamt{\"u}rk, 2003]{A:mip}
Atamt{\"u}rk, A. (2003).
\newblock On the facets of mixed--integer knapsack polyhedron.
\newblock {\em Mathematical Programming}, 98:145--175.

\bibitem[Atamt{\"u}rk et~al., 2016]{AGK:3part-fc}
Atamt{\"u}rk, A., G{\'o}mez, A., and K{\"u}{\c{c}}{\"u}kyavuz, S. (2016).
\newblock Three-partition flow cover inequalities for constant capacity
  fixed-charge network flow problems.
\newblock {\em Networks}, 67:299--315.

\bibitem[Atamt{\"u}rk and G{\"u}nl{\"u}k, 2007]{AG:mixset}
Atamt{\"u}rk, A. and G{\"u}nl{\"u}k, O. (2007).
\newblock Network design arc set with variable upper bounds.
\newblock {\em Networks}, 50:17--28.

\bibitem[Atamt{\"u}rk and G{\"u}nl{\"u}k, 2010]{AG:mingle}
Atamt{\"u}rk, A. and G{\"u}nl{\"u}k, O. (2010).
\newblock Mingling: Mixed-integer rounding with bounds.
\newblock {\em Mathematical Programming}, 123:315--338.

\bibitem[Atamt{\"u}rk and G{\"u}nl{\"u}k, 2018]{AG:capacity}
Atamt{\"u}rk, A. and G{\"u}nl{\"u}k, O. (2018).
\newblock A note on capacity models for network design.
\newblock {\em Operations Research Letters}, 46:414--417.

\bibitem[Atamt{\"u}rk and Kianfar, 2012]{AK:nstepmingle}
Atamt{\"u}rk, A. and Kianfar, K. (2012).
\newblock $n$-step mingling inequalities: new facets for the mixed-integer
  knapsack set.
\newblock {\em Mathematical Programming}, 132:79--98.

\bibitem[Atamt{\"u}rk et~al., 2001]{ANS:avub}
Atamt{\"u}rk, A., Nemhauser, G.~L., and Savelsbergh, M. W.~P. (2001).
\newblock Valid inequalities for problems with additive variable upper bounds.
\newblock {\em Mathematical Programming}, 91:145--162.

\bibitem[Atamt{\"u}rk and Rajan, 2002]{AR:arc}
Atamt{\"u}rk, A. and Rajan, D. (2002).
\newblock On splittable and unsplittable flow capacitated network design
  arc-set polyhedra.
\newblock {\em Mathematical Programming}, 92:315--333.

\bibitem[Avella et~al., 2007]{AMS:metric}
Avella, P., Mattia, S., and Sassano, A. (2007).
\newblock Metric inequalities and the network loading problem.
\newblock {\em Discrete Optimization}, 4:103--114.

\bibitem[Balakrishnan et~al., 1991]{BMSW:nexp}
Balakrishnan, A., Magnanti, T.~L., Shulman, A., and Wong, R.~T. (1991).
\newblock Models for planning capacity expansion in local access
  telecommunication networks.
\newblock {\em Annals of Operations Research}, 33:239--284.

\bibitem[Balakrishnan et~al., 1995]{BMW:local-access}
Balakrishnan, A., Magnanti, T.~L., and Wong, R.~T. (1995).
\newblock A decomposition algorithm for local access telecommunication network
  expansion planning.
\newblock {\em Operations Research}, 43:58--76.

\bibitem[Barahona, 1996]{B:nd-cut}
Barahona, F. (1996).
\newblock Network design using cut inequalities.
\newblock {\em SIAM Journal on Optimization}, 6:823--837.

\bibitem[Barnhart et~al., 2000]{BHV:intmulticom}
Barnhart, C., Hane, C.~A., and Vance, P.~H. (2000).
\newblock Using branch-and-price-and-cut to solve origin-destination integer
  multicommodity flow problems.
\newblock {\em Operations Research}, 48:318--326.

\bibitem[Bienstock et~al., 1998]{BCGT:cap-inst}
Bienstock, D., Chopra, S., G{\"u}nl{\"u}k, O., and Tsai, C.-Y. (1998).
\newblock Minimum cost capacity installation for multicommodity networks.
\newblock {\em Mathematical Programming}, 81:177--199.

\bibitem[Bienstock and G{\"u}nl{\"u}k, 1996]{BG:cnd}
Bienstock, D. and G{\"u}nl{\"u}k, O. (1996).
\newblock Capacitated network design - {P}olyhedral structure and computation.
\newblock {\em INFORMS Journal on Computing}, 8:243--259.

\bibitem[Brockm{\"u}ller et~al., 2004]{BGW:private}
Brockm{\"u}ller, B., G{\"u}nl{\"u}k, O., and Wolsey, L.~A. (2004).
\newblock Designing private line networks -- {P}olyhedral analysis and
  computation.
\newblock {\em Transactions on Operational Research}, 16:7--24.

\bibitem[Chopra et~al., 1998]{CGS:nd-batch}
Chopra, S., Gilboa, I., and Sastry, S.~T. (1998).
\newblock Source sink flows with capacity installation in batches.
\newblock {\em Discrete Applied Mathematics}, 85:165--192.

\bibitem[Costa et~al., 2009]{CCG:mc-netdes}
Costa, A.~M., Cordeau, J.-F., and Gendron, B. (2009).
\newblock Benders, metric and cutset inequalities for multicommodity
  capacitated network design.
\newblock {\em Computational Optimization and Applications}, 42:371--392.

\bibitem[Dahl and Stoer, 1998]{DahlStoer1998}
Dahl, G. and Stoer, M. (1998).
\newblock A cutting plane algorithm for multicommodity survivable network
  design problems.
\newblock {\em INFORMS Journal on Computing}, 10:1--11.

\bibitem[Dash and G{\"u}nl{\"u}k, 2006]{DG:twostep}
Dash, S. and G{\"u}nl{\"u}k, O. (2006).
\newblock Valid inequalities based on simple mixed-integer sets.
\newblock {\em Mathematical Programming}, 105:29--53.

\bibitem[Dash et~al., 2016]{DGW:cont-knap}
Dash, S., G{\"u}nl{\"u}k, O., and Wolsey, L.~A. (2016).
\newblock The continuous knapsack set.
\newblock {\em Mathematical Programming}, 155:471--496.

\bibitem[Davarnia et~al., 2019]{JP:train}
Davarnia, D., Richard, J.-P.~P., \.{I}{\c{c}}y{\"u}z Ay, E., and Taslimi, B.
  (2019).
\newblock Network models with unsplittable node flows with application to unit
  train scheduling.
\newblock {\em Operations research}, 67:1053--1068.

\bibitem[Florian et~al., 1976]{FBFGN:train}
Florian, M., Bushell, G., Ferland, J., Gu{\'e}rin, G., and Nastansky, L.
  (1976).
\newblock The engine scheduling problems in a railway network.
\newblock {\em INFOR}, 14:121--138.

\bibitem[Gavish and Altinkemer, 1990]{GA:backbone}
Gavish, B. and Altinkemer, K. (1990).
\newblock Backbone network design tools with economic tradeoffs.
\newblock {\em ORSA Journal on Computing}, 2:58--76.

\bibitem[G{\"u}nl{\"u}k, 1999]{G:bc-nd}
G{\"u}nl{\"u}k, O. (1999).
\newblock A branch-and-cut algorithm for capacitated network design problems.
\newblock {\em Mathematical Programming}, 86:17--39.

\bibitem[G{\"{u}}nl{\"{u}}k, 2007]{DBLP:journals/siamdm/Gunluk07}
G{\"{u}}nl{\"{u}}k, O. (2007).
\newblock A new min-cut max-flow ratio for multicommodity flows.
\newblock {\em {SIAM} J. Discrete Mathematics}, 21:1--15.

\bibitem[Hamid and Agarwal, 2015]{A:Hamid:2015}
Hamid, F. and Agarwal, Y.~K. (2015).
\newblock Solving the two-facility network design problem with 3-partition
  facets.
\newblock {\em Networks}, 66:11--32.

\bibitem[{\SortNoop{Hoesel}van Hoesel} et~al., 2002]{HKLS:edge-cap}
{\SortNoop{Hoesel}van Hoesel}, S. P.~M., Koster, A. M. C.~A., {van de Leensel},
  R. L. M.~J., and Savelsbergh, M. W.~P. (2002).
\newblock Polyhedral results for the edge capacity polytope.
\newblock {\em Mathematical Programming}, 92:335--358.

\bibitem[Iri, 1971]{I:metric}
Iri, M. (1971).
\newblock On an extension of the max-flow min-cut theorem to multicommodity
  flows.
\newblock {\em Journal of the Operations Research Society of Japan},
  13:129--135.

\bibitem[Labb{\'e} and Yaman, 2004]{LY:project}
Labb{\'e}, M. and Yaman, H. (2004).
\newblock Projecting the flow variables for hub location problems.
\newblock {\em Networks}, 44:84--93.

\bibitem[Magnanti and Mirchandani, 1993]{MM:st-path}
Magnanti, T.~L. and Mirchandani, P. (1993).
\newblock Shortest paths, single origin-destination network design, and
  associated polyhedra.
\newblock {\em Networks}, 23:103--121.

\bibitem[Magnanti et~al., 1993]{MMV:conv-2core}
Magnanti, T.~L., Mirchandani, P., and Vachani, R. (1993).
\newblock The convex hull of two core capacitated network design problems.
\newblock {\em Mathematical Programming}, 60:233--250.

\bibitem[Magnanti et~al., 1995]{MMV:modeling-2fnlp}
Magnanti, T.~L., Mirchandani, P., and Vachani, R. (1995).
\newblock Modeling and solving the two--facility capacitated network loading
  problem.
\newblock {\em Operations Research}, 43:142--157.

\bibitem[Magnanti and Wong, 1984]{MW:nd}
Magnanti, T.~L. and Wong, R.~T. (1984).
\newblock Network design and transportation planning: Models and algorithms.
\newblock {\em Transportation Science}, 18:1--55.

\bibitem[Marchand and Wolsey, 2001]{MW:MIR}
Marchand, H. and Wolsey, L.~A. (2001).
\newblock Aggregation and mixed integer rounding to solve {MIPs}.
\newblock {\em Operations Research}, 49:363--371.

\bibitem[Mattia, 2012]{M:sep-metric}
Mattia, S. (2012).
\newblock Separating tight metric inequalities by bilevel programming.
\newblock {\em Operations Research Letters}, 40:568 -- 572.

\bibitem[Minoux, 1989]{M:netsyn}
Minoux, M. (1989).
\newblock Network synthesis and optimum network design problems: Models,
  solution methods and applications.
\newblock {\em Networks}, 19:313--360.

\bibitem[Mirchandani, 2000]{M:proj-netload}
Mirchandani, P. (2000).
\newblock Projections of the capacitated network loading problem.
\newblock {\em European Journal of Operational Research}, 122:534--560.

\bibitem[Nemhauser and Wolsey, 1988]{NW:IPbook}
Nemhauser, G.~L. and Wolsey, L.~A. (1988).
\newblock {\em Integer and Combinatorial Optimization}.
\newblock John Wiley and Sons, New York.

\bibitem[Onaga and Kakusho, 1971]{OK:metric}
Onaga, K. and Kakusho, O. (1971).
\newblock On feasibility conditions of multi-commodity flows in networks.
\newblock {\em Transactions on Circuit Theory}, 18:425--429.

\bibitem[Ozbaygin et~al., 2018]{OKY:splitdelivery}
Ozbaygin, G., Karasan, O., and Yaman, H. (2018).
\newblock New exact solution approaches for the split delivery vehicle routing
  problem.
\newblock {\em EURO Journal on Computational Optimization}, 6:85--115.

\bibitem[Pochet and Wolsey, 1995]{PW:div}
Pochet, Y. and Wolsey, L.~A. (1995).
\newblock Integer knapsack and flow covers with divisible coefficients:
  Polyhedra, optimization, and separation.
\newblock {\em Discrete Applied Mathematics}, 59:57--74.

\bibitem[Raack et~al., 2011]{RKOW:cutset}
Raack, C., Koster, A.~M., Orlowski, S., and Wess{\"a}ly, R. (2011).
\newblock On cut-based inequalities for capacitated network design polyhedra.
\newblock {\em Networks}, 57:141--156.

\bibitem[Wolsey, 1998]{W:IPbook}
Wolsey, L.~A. (1998).
\newblock {\em Integer Programming}.
\newblock John Wiley and Sons, New York.

\bibitem[Wolsey and Yaman, 2016]{WY:contknap-div}
Wolsey, L.~A. and Yaman, H. (2016).
\newblock Continuous knapsack sets with divisible capacities.
\newblock {\em Mathematical Programming}, 156:1--20.

\bibitem[Yaman, 2007]{Y:knap-cover}
Yaman, H. (2007).
\newblock The integer knapsack cover polyhedron.
\newblock {\em SIAM Journal on Discrete Mathematics}, 21:551--572.

\bibitem[Yaman, 2013]{Y:arcset}
Yaman, H. (2013).
\newblock The splittable flow arc set with capacity and minimum load
  constraints.
\newblock {\em Operations Research Letters}, 41:556--558.

\end{thebibliography}

\end{document}